\documentclass[11pt,reqno]{amsart}
\usepackage{amsfonts}
\usepackage{fancyhdr}
\usepackage{times}
\usepackage[T1]{fontenc}
\usepackage{mathrsfs}
\usepackage{latexsym}
\usepackage[dvips]{graphics}
\usepackage{epsfig}
\usepackage{hyperref, amsmath, amsthm, amsfonts, amscd, flafter,epsf}
\usepackage{amsmath,amsfonts,amsthm,amssymb,amscd}
\input amssym.def
\input amssym.tex
\usepackage{color}
\usepackage[parfill]{parskip}
 \usepackage{mathrsfs} 
\usepackage{geometry}                
\usepackage{epstopdf}

\usepackage{verbatim}

\newtheorem{thm}{Theorem}
\newtheorem{cor}[thm]{Corollary}
\newtheorem{con}[thm]{Conjecture}
\newtheorem{lem}[thm]{Lemma}
\newtheorem*{thm*}{Theorem}
\newtheorem*{con*}{Conjecture}
\newtheorem*{lem*}{Lemma}

\newtheorem{defn}{Definition}[section]
\newtheorem{rem*}{Remark}

\begin{document}
 \baselineskip=17pt
\hbox{}
\medskip

\title
[index and Dedekind sum] {note on the index conjecture in zero-sum
theory
\\and its connection to a Dedekind-type sum}

\author{Fan Ge}

\email{fan.ge@rochester.edu}

\address{Department of Mathematics, University of Rochester, Rochester, NY, United States}

\thanks{This work was partially supported by NSF grant DMS-1200582.}

\maketitle

\begin{abstract}
Let $S=(a_1)\cdots(a_k)$ be a minimal zero-sum sequence over a
finite cyclic group $G$. The index conjecture states that if $k=4$
and $\gcd(|G|,6)=1$, then $S$ has index 1. In this note we study the
index conjecture and connect it to a Dedekind-type sum. In
particular we reprove a special case of the conjecture when $|G|$ is
prime.
\end{abstract}

\section{Introduction}

Throughout this paper let $G$ be a finite cyclic group of order $n$,
written additively.  By a \emph{sequence $S$ of length $k$ over $G$}
we mean a sequence with $k$ elements, each of which is in $G$. We
write $(a_1)\cdots(a_k)$ for such a sequence, rather than $a_1$,
\dots, $a_k$. A sequence $S$ is said to be a \emph{zero-sum}
sequence if $\sum_i a_i=0$. It is a \emph{minimal} zero-sum sequence
if it is a zero-sum sequence but no proper nontrivial subsequence of
it is. Given any $g$ a generator of $G$, we can write
$S=(x_1g)\cdots(x_kg)$ for some natural numbers $x_1, \dots, x_k$,
where by $x_ig$ we mean the sum $g+g+ \cdots +g$ with $x_i$ terms.
We also require the following definition.

\begin{defn}\textnormal{For a sequence over $G$
$$S=(x_1g)\cdots(x_kg), \qquad \textnormal{where} \ 1\le x_1,...,x_k\le n\,,$$
define the $g$-norm of $S$ to be
$\|S\|_g=\frac{\sum_{i=1}^{k}x_i}{n}$. The \emph{index} of $S$ is
defined by
$$\textnormal{ind}(S)=\min\|S\|_g,$$ where the minimum is taken
over all generators $g$ of $G$.}
\end{defn}

The index of a sequence is a crucial invariant in the theory of
zero-sum sequences over cyclic groups. It was  introduced by
Kleitman and Lemke~\cite{kl} and then became an important tool in the
study of zero-sum sequences and related topics (see, for example, Geroldinger~\cite{g90}  and Gao~\cite{gao}).
An important question is to determine the pairs
$(k, n)$ for which every minimal zero-sum sequence $S$ of length $k$
over $G$ has index 1. The cases  $k\ne 4$ or
$\gcd(n, 6)\ne 1$ have been settled (see~\cite{G}, \cite{P}, \cite{sc}, \cite{Y}).
Hence, the only remaining case   is when   $k=4$ and, at the same time,  $\gcd(n,
6)=1$. The following conjecture is widely held.

\begin{con}\label{conj}
Let $G$ be a finite cyclic group such that $\gcd(|G|, 6) = 1$. Then
every minimal zero-sum sequence $S$ over $G$ of length 4 has
\textnormal{ind}($S$) = 1.
\end{con}

\begin{rem*}
\textnormal{One can show that, for a minimal zero-sum sequence $S$
of length 4 over $G$, we have either $\textnormal{ind}(S)=1$ or
$\textnormal{ind}(S)=2$, and moreover, $\textnormal{ind}(S)=2$ if
and only if $\|S\|_g=2$ for all generators $g$ of $G$. This is
because for such $S$ we have by definition that $\|S\|_g$ is equal
to 1, 2, or 3; but if $\|S\|_g=3$ for some generator $g$, then it is
easy to see that $\|S\|_{-g}=1$, where $-g$ is also a generator.}
\end{rem*}

In~\cite{LPYZ}, Y. Li et al. proved that Conjecture~\ref{conj} is
true for when $n$ is a prime power. Later it was proved for the case
when $n$ has two distinct prime factors (see~\cite{LP}
and~\cite{XS}). We also have the following contribution of Shen, Xia and Li~\cite{SXL}.

\begin{thm}\label{X's thm}
Let $G = \langle g\rangle$ be a finite cyclic group of order $|G| = n$ such that $\gcd(n, 6) =
1$. Let $S = (x_1g)(x_2g)(x_3g)(x_4g)$ be a minimal zero-sum
sequence over $G$, where  $\gcd(n, x_1, x_2, x_3, x_4) = 1$, and for
some $i\in [1, 4]$ we have $\gcd(n, x_i)> 1$. Then
\textnormal{ind}($S$) = 1.
\end{thm}

The purpose of this paper is to give a new proof of Conjecture~\ref{conj}
in the case when $n$ is a prime and at the same time to point out a connection between
the conjecture and the theory of Dedekind sums.

Recall that the classical Dedekind sum $s(h, k)$, where $k\ge 1$ and
$\gcd(h,k)=1$, is defined as
$$
s(h, k)=\sum_{r=1}^k \frac{r}{k}\Big(\frac{hr}{k}-\Big\lfloor
\frac{hr}{k}\Big\rfloor - \frac{1}{2}\Big).
$$
This has various natural generalizations (for example, see ~\cite{Z}
and ~\cite{RG}). Here we shall consider another
Dedekind-type sum. For $k\ge 1$ and $\gcd(h,k)=1$ let
$$t(h, k)=\sum_{\substack{1\le r\le k\\(r,k)=1}}
\frac{r}{k}\Big(\frac{hr}{k}-\Big\lfloor \frac{hr}{k}\Big\rfloor -
\frac{1}{2}\Big).$$
Clearly $t(h,k)$ can be expressed in terms of classical Dedekind
sums using M\"{o}bius convolution.
In this note we  connect the index conjecture with $t(h,k)$.

\begin{thm}\label{thm i}
\textnormal{(i)} Suppose that $n$ is the smallest integer for which
Conjecture~\ref{conj} fails. Let $S=(x_1)(x_2)(x_3)(x_4)$ be
a minimal zero-sum sequence over $G\cong \mathbb{Z}/n$ with \textnormal{ind}$(S)=2$. Then we have $\gcd(n, x_i)=1$ for all $i$.

\textnormal{(ii)} Let $n$ be a positive integer for which Conjecture~\ref{conj}
fails, and let $S=(x_1)(x_2)(x_3)(x_4)$ be a minimal zero-sum
sequence over $\mathbb{Z}/n$ with \textnormal{ind}$(S)=2$. If
$\gcd(n,x_i)=1$ for all $i$, then we have
\begin{align}\label{eq t(h,k) in thm}
t(\overline{x_{\sigma(1)}}\
x_{\sigma(2)},n)=t(\overline{x_{\sigma(3)}}\ x_{\sigma(4)},n)
\end{align}
for any permutation $\sigma$ of $\{1, 2, 3, 4\}$, where
$\overline{x}$ denotes the inverse of $x$ in the multiplicative
group $(\mathbb{Z}/n)^*$.\end{thm}

\begin{cor}\label{cor} Conjecture \ref{conj} is true if $G$ has prime order.\end{cor}

\section{proof of theorem \ref{thm i} and corollary \ref{cor}}

For integers $x$ and $y>0$, let $(x)_y$ denote the least nonnegative
residue of $x$ mod $y$. For $z\in \mathbb Z/y$, we may view $z$ as
an integer and define $(z)_y$ similarly.

\begin{lem} \label{lem coeff} Let $S=(x_1)(x_2)(x_3)(x_4)$ be a sequence over $G=\mathbb{Z}/n$.
Given any generator $g$ in $G$, write $S=(y_1g)(y_2g)(y_3g)(y_4g)$
for $1\le y_i \le n$. Then we have $y_i=(g^{-1}x_i)_n$ for $i=1,
..., 4$, where $g^{-1}$ is the inverse of $g$ in the multiplicative
group $(\mathbb Z/n)^*$.
\end{lem}

\proof For any $i=1, ..., 4$, we have $x_i=y_i g$ in $G$. Hence
$x_ig^{-1}=y_i$ in $G$. Therefore, $y_i=(g^{-1}x_i)_n$.

\emph{Proof of Theorem \ref{thm i}:}
(i)
Suppose that $n$ is the
smallest positive integer for which Conjecture \ref{conj} fails, and
$S=(x_1)(x_2)(x_3)(x_4)$ is a minimal zero-sum sequence over
$\mathbb{Z}/n$ with ind$(S)$=2. In view of Theorem~\ref{X's thm} we
know that either $\gcd(n, x_1, x_2, x_3, x_4) >1$, or $\gcd(n,
x_i)=1$ for all $i$. We will show that the first case does not
happen.

Suppose that $\gcd(n, x_1, x_2, x_3, x_4) =d >1$. Let $m=n/d$, and
consider
$S'=(\tfrac{x_1}{d})(\tfrac{x_2}{d})(\tfrac{x_3}{d})(\tfrac{x_4}{d})$
over $\mathbb{Z}/m$. It is easy to check that $S'$ is a minimal
zero-sum sequence. We will show that ind$(S')=2$, so that
Conjecture~\ref{conj} fails for $m$. This contradicts the assumption
that $n=md$ is the smallest integer for which Conjecture~\ref{conj}
fails.

By Remark 1, to prove ind$(S')=2$ it suffices to show that
for every $h\in (\mathbb{Z}/m)^*$, we have $\|S'\|_h=2$. By Lemma
\ref{lem coeff} and the definition of $\|S'\|_h$, this is nothing but $\sum_i
(h^{-1}\frac{x_i}{d})_{m}=2m$. Note that
$(h^{-1}\frac{x_i}{d})_{m}=(h^{-1}x_i)_n/d$, so it suffices to show
that $\sum_i (h^{-1}x_i)_n/d=2m$ or, equivalently, that
$\sum_i(h^{-1}x_i)_n=2n$. But this last assertion is guaranteed by Remark 1 and the assumption that ind$(S)=2$.

It  now follows that we   have  $\gcd(n, x_1, x_2, x_3, x_4) = 1$, and thus
 $\gcd(n, x_i)=1$ for all $i$ by Theorem \ref{X's thm}. This proves part (i).

\medskip

(ii) Consider the sum
\begin{align} \label{wt sum}
W := \sum_{\substack{0 < g < n \\(g,n)=1}}
\bigg(\sum_{i=1}^4(gx_i)_n\bigg)\cdot \bigg(\sum_{i=1}^4
(-1)^i(gx_i)_n\bigg).
\end{align}
By multiplying out the product, we see that
\begin{align*}
W = 2 \sum_{\substack{0< g < n \\(g,n)=1}}  \bigg((gx_2)_n(gx_4)_n
-(gx_1)_n(gx_3)_n \bigg) +  \sum_{\substack{0< g < n \\(g,n)=1}}
\bigg( (gx_2)_n^2 + (gx_4)_n^2 - (gx_1)_n^2 - (gx_3)_n^2\bigg).
\end{align*}

Recall that $\gcd(n, x_i)=1$ for all $i\in [1, 4]$. Therefore, for
any $i, j\in [1, 4]$ and any integer $k$,
\begin{align}\label{eq power sum}
\sum_{\substack{0< g < n \\(g,n)=1}}(gx_i)_n^k = \sum_{\substack{0<
g < n \\(g,n)=1}} (gx_j)_n^k.
\end{align}
Hence we see that
\begin{align}\label{w cross term}
W = 2 \sum_{\substack{0< g < n \\(g,n)=1}}  \big((gx_2)_n(gx_4)_n
-(gx_1)_n(gx_3)_n \big).
\end{align}

On the other hand, since ind$(S)=2$, by Remark 1 and Lemma
\ref{lem coeff} we have $\sum_{i=1}^4(gx_i)_n=2n$ for all $g$
coprime to $n$. It then follows from equation~(\ref{wt sum}) that
\begin{align*}
W & = \sum_{\substack{0< g < n \\ (g,n)=1}} 2n \bigg(\sum_{i=1}^4 (-1)^i(gx_i)_n\bigg)
\\
& = 2n \bigg( \sum_{\substack{0< g < n \\ (g,n)=1}} \big( (gx_2)_n +
(gx_4)_n - (gx_1)_n - (gx_3)_n \big) \bigg),
\end{align*}
and, in view of~(\ref{eq power sum}), this is
$$
W = 0.
$$
Combining this with (\ref{w cross term}) we obtain that
\begin{align}\label{necessary condition}
\sum_{\substack{0< g<n\\(g,n)=1}}
(gx_2)_n(gx_4)_n=\sum_{\substack{0< g<n\\(g,n)=1}}
(gx_1)_n(gx_3)_n,\end{align}

or equivalently,
\begin{align}\label{necessary condition 2}
\sum_{\substack{0< g<n\\(g,n)=1}}
g(g\overline{x_2}x_4)_n=\sum_{\substack{0< g<n\\(g,n)=1}}
g(g\overline{x_1}x_3)_n.\end{align}

\medskip

Recall that
\begin{align*}
t(h, k) & =\sum_{\substack{0<r<k \\(r,k)=1}} \frac{r}{k}\Big(\frac{hr}{k}-\Big\lfloor
\frac{hr}{k}\Big\rfloor - \frac{1}{2}\Big)
\\
& = \sum_{\substack{0< r<k\\(r,k)=1}} \frac{r}{k}\Big(\frac{(hr)_k}{k}- \frac{1}{2}\Big)
\\
& = \frac{1}{k^2}\sum_{\substack{0<r<k\\(r,k)=1}} r(hr)_k - \frac{1}{2k}\sum_{\substack{0< r<k\\(r,k)=1}} r.
\end{align*}

Hence, by ~(\ref{necessary condition 2}) we see that
$t(\overline{x_2}x_4,n)=t(\overline{x_1}x_3,n)$, and
equation~(\ref{eq t(h,k) in thm}) follows by the symmetry of
$x_i$'s. This completes the proof of Theorem~\ref{thm i}. \qed

\medskip

To prove Corollary \ref{cor} we require the following lemma on Dedekind sums.
\begin{lem}
 Let $k$ be a prime. If $s(h_1, k)=s(h_2, k)$, then we have $h_1 \equiv h_2 \pmod
k$ or $h_1 \equiv \overline{h_2} \pmod k$.
\end{lem}
\proof
This is Corollary 1.2 in~\cite{JRW2}, but
we sketch the proof for completeness.
The equality $s(h_1, k)=s(h_2, k)$ implies that $12ks(h_1, k)\equiv
12ks(h_2, k) \pmod{k}$. From the well-known reciprocity law for
Dedekind sums (see Chapter 3 of~\cite{A}, for example), it is not hard to derive that
$12ks(h,k)\equiv h+\overline{h} \pmod{k}.$
The result then follows by a straightforward computation.

\medskip

\emph{Proof of Corollary \ref{cor}:} Let $n=p$ be a
prime, and let $S=(x_1)\cdots(x_4)$ be a minimal zero-sum sequence
over $\mathbb{Z}/p$. Suppose that ind$(S)=2$. By Theorem~\ref{thm
i} we have $t(\overline{x_2}x_4, p)=t(\overline{x_1}x_3, p)$, or
equivalently, $s(\overline{x_2}x_4, p)=s(\overline{x_1}x_3, p)$.

Since $s(\overline{x_2}x_4, p)=s(\overline{x_1}x_3, p)$, by the above lemma we
see that
$$
\overline{x_2}x_4 \equiv \overline{x_1}x_3 \pmod{p} \ \ \textrm{ or }\ \
\overline{x_2}x_4 \equiv \overline{x_3}x_1 \pmod{p}.
$$
Without loss of generality, suppose that $\overline{x_2}x_4 \equiv
\overline{x_1}x_3 \pmod{p}.$ Then we have $x_1x_4 \equiv x_2x_3
\pmod{p}$. It follows that
$$
0\equiv
x_1x_4-x_2x_3\equiv x_1x_4+x_2(x_1+x_2+x_4)=(x_2+x_1)(x_2+x_4)
\pmod{p}.
$$
Hence, we have $p\mid x_2+x_1$ or $p\mid x_2+x_4$, both of which
contradicts the fact that $S$ is a minimal zero-sum sequence over
$\mathbb{Z}/p$. This completes our proof of Corollary \ref{cor}.
\qed

\medskip

\section{further discussion}
An interesting and  natural  question about Dedekind sums is:  for
which $h_1, h_2$ and $k$  do we have
$$
s(h_1,k)=s(h_2,k)\ ?
$$
(For example, this question arises in connection   with the Heegaard
Floer Homology ~\cite{JRW}.) A similar question is to determine when
we have $t(h_1,k)=t(h_2,k)$. This is also interesting in view of
Theorem~\ref{thm i}. Unfortunately, for now we are unable to answer
either of these two questions. (See also~\cite{JRW2}, \cite{Gir}
and~\cite{Tsu}.)

We end our note by the following remark. According to
Theorem~\ref{thm i}, a possible way to prove the index conjecture
is to show that there exists at least one permutation $\sigma$ of $\{1, 2, 3, 4\}$
for which equation~(\ref{eq t(h,k) in thm}) does not hold.
Note that for a given $n$ and minimal zero-sum sequence $S=(x_1)(x_2)( x_3)( x_4)$ it is not necessarily the case that ~(\ref{eq t(h,k) in thm}) will hold for \emph{no} permutation of $\{1,2,3,4\}$.
For example, if $n=25$, $x_1=18$, $x_2=4$, $x_3=2$
and $x_4=1$, then we have ind$(S)=1$, so  Conjecture 1 is true. However,
$t(\overline{x_2}x_4,n)=t(\overline{x_1}x_3,n)$, but
$t(\overline{x_2}x_1,n)\ne t(\overline{x_4}x_3,n)$.

\medskip

\section*{Acknowledgement}

The author would like to thank Professors Noam Elkies, Alfred
Geroldinger, Steve Gonek and Gerry Myerson for helpful
communications. Special thanks to Professor Steve Gonek for
improving the presentation of the paper.

\end{document}